\documentclass[a4paper,12pt]{amsart}

\usepackage{amsmath,times,epsfig,amssymb,amsbsy,amscd,amsfonts,amstext,color,bm}

\usepackage[all]{xy}
\usepackage{bm}
\usepackage{comment}
\usepackage{enumerate}
\usepackage{circledsteps}
\usepackage{geometry}

\usepackage[dvipdfmx,svgnames]{xcolor}
\usepackage{tikz}
\tikzset{v/.style={circle, draw, inner sep=2pt, minimum size=6pt, fill=white}}

\theoremstyle{plain}
\newtheorem{theorem}{Theorem}[section]
\newtheorem*{theorem*}{Theorem}
\newtheorem*{theoremA*}{Theorem A}
\newtheorem*{theoremB*}{Theorem B}
\newtheorem{corollary}[theorem]{Corollary}

\newtheorem{conjecture}[theorem]{Conjecture}

\theoremstyle{definition}
\newtheorem{definition}[theorem]{Definition}
\newtheorem{remark}[theorem]{Remark}

\newtheorem{proposition}[theorem]{Proposition}

\newcommand{\Hom}{\operatorname{Hom}}

\newcommand{\Tor}{\operatorname{Tor}}

\newcommand{\codim}{\operatorname{codim}}

\newcommand{\bC}{\mathbb{C}}

\newcommand{\bZ}{\mathbb{Z}}

\newcommand{\cA}{\mathcal{A}}
\newcommand{\cB}{\mathcal{B}}

\def\qed{\hfill $\Box$}

\title{Even torsions in the homology group of the Milnor fiber boundary of hyperplane arrangements in $\mathbb{C}^3$}
\author{Sakumi Sugawara}
\date{\today}
\address{Department of Mathematics, Faculty of Science, Hokkaido University, North 10, West 8, Kita-ku, Sapporo 060-0810, JAPAN. }
\email{sugawaras@math.sci.hokudai.ac.jp}

\begin{document}
\maketitle

\begin{abstract}
We study the homology group of the Milnor fiber boundary of a hyperplane arrangement in $\bC^{3}$.
By the work of N\'emethi--Szil\'ard, the homeomorphism type of the Milnor fiber boundary is combinatorially determined, and an explicit formula for the first Betti number is known.
However, the torsion part of the first homology group is poorly understood.
In this paper, under some conditions, we prove that the number of even-order torsion summands of the first homology group is greater than or equal to the Euler characteristic of the projectivized complement.
\end{abstract}


\section{Introduction}
A finite set $\cA = \{H_1, \ldots, H_{n}\}$ of linear hyperplanes in a complex Euclidean space is called a hyperplane arrangement.
Hyperplane arrangements are studied from various viewpoints, including algebraic geometry, combinatorics, and topology \cite{orl-ter}.
One of the major problems in the topology of hyperplane arrangements is whether several topological invariants are combinatorially determined (that is, by the information of the intersection poset).
It is known that the cohomology ring of the complement is combinatorially determined \cite{orl-sol}.
There are a lot of studies on a combinatorial description of the Betti numbers of the Milnor fiber, but a complete solution has not yet been achieved (see \cite{pap-suc-17} for the recent progress on this topic).

In this paper, we focus on the Milnor fiber boundary of arrangements in $\bC^3$. The Milnor fiber boundary is a counterpart of the singularity link for an isolated hypersurface singularity. 
N\'emethi--Szli\'ard studied the Milnor fiber boundary of a non-isolated surface singularity in \cite{nem-szi}. They proved that the Milnor fiber boundary is a plumbed manifold (see also \cite{mic-pic}) and gave an algorithm to obtain the plumbing graph.
By applying their result to hyperplane arrangements, it follows that the homeomorphism type of the Milnor fiber boundary is combinatorially determined.
N\'emethi--Szil\'ard also gave the following explicit combinatorial description for the first Betti number (see also Theorem 7.8 in \cite{suc-mil}).
\begin{theorem}\label{thm:mfb_betti} (Theorem 19.10.2 in \cite{nem-szi})
Let $\cA$ be a hyperplane arrangement with $n$ hyperplanes. Then, the first Betti number of the Milnor fiber boundary $\partial F$ is described as
\[
b_{1} (\partial F) = \sum_{X \in L_{2} (\cA)} (1 + (|\cA_{X}|-2) \gcd(|\cA_{X}|, n) ).
\]
(The definitions of notations will be given in Section \ref{sec:arr}.)
\end{theorem}

However, the first homology group of the Milnor fiber boundary has a non-trivial torsion in general. Providing an explicit formula describing the torsion part has also been posed as a problem by N\'emethi--Szil\'ard (see Problem 24.4.19 in \cite{nem-szi}). 
Although there are only a few studies on this topic, the following results and conjectures are known.

\begin{theorem}\label{thm:gen} (Theorem 1.3 in \cite{sug-mfb})
Let $\cA$ be a generic arrangement with $n$ hyperplanes in $\bC^{3}$.
Then, the first homology group of the Milnor fiber boundary is described as
\[
H_{1} (\partial F; \bZ) \cong \bZ^{n(n-1)/2} \oplus \bZ_{n}^{(n-2)(n-3)/2}.
\]
\end{theorem}

\begin{conjecture} (Conjecture 1.4 in \cite{sug-mfb})
Suppose that $(|\cA_{X}|-2)( \gcd(|\cA_X| , n)- 1) =0$ for each $X \in L_{2} (\cA)$. Then, the torsion part $\Tor (H_{1} (\partial F); \bZ)$ satisfies
\[
\Tor (H_{1} (\partial F); \bZ) \cong (\bZ_{n})^{\oplus \chi (U)},
\]
where $\chi (U)$ is the Euler characteristic of the projectivized complement.
\end{conjecture}

Building on these observations, this paper provides a lower bound for the number of even torsion summands for a somewhat larger class of arrangements than generic ones.
The following is the main result of this paper.

\begin{theorem} \label{thm:main}
Let $\cA$ be a central hyperplane arrangement whose number of hyperplanes is even.
Suppose that $(|\cA_X|-2)(\gcd(|\cA_X|, n) -1) =0$ for each $X \in L_2 (\cA)$.
Then, 
\[
\dim_{\bZ_2} (\Tor(H_1 (\partial F; \bZ)) \otimes \bZ_2 ) \geq \chi (U).
\]
\end{theorem}
Note that when $n$ is even, the condition $(|\cA_X|-2)(\gcd(|\cA_X|, n) -1) =0$ for each $X \in L_2 (\cA)$ is equivalent to the condition ``if $|\cA_X|$ is not $2$, then $\cA_X$ is odd''.
When $|\cA_X|=2$ for every $X \in L_{2} (\cA)$ (that is, when the arrangement is generic), the result agrees with Theorem \ref{thm:gen}.
In particular, the inequality becomes an equality.

The proof uses techniques for analyzing the $2$-torsion in the homology groups of double covers developed in \cite{yos-double, isy}, together with results on the resonance varieties of boundary manifolds studied in \cite{coh-suc-proj, coh-suc-bound}.


This paper is organized as follows. In Section 2, we review the topology of hyperplane arrangements, with particular emphasis on boundary manifolds and Milnor fiber boundaries. Section 3 summarizes the results and techniques for studying the homology groups of the tower of double coverings. The proof of the main theorem is given in Section 4. In Section 5, we present a list of computations of the homology groups of Milnor fiber boundaries for the case of eight hyperplanes.
\vspace{3mm}
\\
\textbf{Acknowledgements}. The author would like to thank Masahiko Yoshinaga for the helpful discussion and Yongqiang Liu for valuable comments on this paper. The author is supported by JSPS KAKENHI 25K23326.

\section{Topology of hyperplane arrangements} \label{sec:arr}
In this section, we summarize fundamental properties of the topology of hyperplane arrangements that will be needed. For details, see \cite{orl-ter, suc-mil, dim-hyp} for example.

\subsection{The characterisitic polynomial}

Let $\cA = \{H_1, \ldots, H_{n}\}$ be a central hyperplane arrangement in $\bC^{\ell+1} $ with the complement $M(\cA) = \bC^{\ell+1} \setminus \bigcup_{i=1}^{n} H_{i}$.
The set of intersections of hyperplanes defines the intersection poset $L(\cA) = \{\bigcap_{H \in \cB} H \mid \cB \subset \cA \}$ with the order $X \leq Y \Leftrightarrow X \supset Y$. Let $L_{k} (\cA) = \{X \in L(\cA) \mid \codim X =k \}$.
For $X \in L(\cA)$, the localization $\cA_{X}$ is defined as a subarrangement $\{H \in \cA \mid H \supset X\}$.

The M\"obius function $\mu: L(\cA) \rightarrow \bZ$ is defined inductively as follows:
\begin{eqnarray*}
\mu (X) = 
\left \{
\begin{array}{ll}
1 & (X= \bC^{\ell+1}), \\
-\sum_{Y < X} \mu (Y) & (X > \bC^{\ell + 1}).
\end{array}
\right .
\end{eqnarray*}
For example, $\mu(X) = -1$ for $X \in L_{1} (\cA)$ and $\mu(X) = |\cA_{X}| -1$ for $X \in L_2 (\cA)$.
The characteristic polynomial $\chi (\cA,t)$ is defined as 
\[
\chi (\cA, t) = \sum_{X \in L(\cA)} \mu(X) t^{\dim X}.
\]
It is well-known that the characteristic polynomial describes the Betti numbers of the complement. 
\begin{proposition}
\[
\sum_{k \geq 0} b_{k} (M)t^{k} = (-t)^{\ell + 1} \chi (\cA, - \frac{1}{t}).
\] 
In particular, $b_{k} (M) = (-1)^{k}\sum_{X \in L_{k}(\cA)} \mu (X)$.
\end{proposition}

\subsection{Milnor fiber boundaries and boundary manifolds}

We fix the defining linear function $\alpha_{i}$ of each hyperplane $H_{i}$.
The product $Q= \prod_{i=1}^{n} \alpha_{i} $ defines a polynomial function $Q: \bC^{\ell+1} \rightarrow \bC$.
Its restriction to the complement $Q: M(\cA) \rightarrow \bC^{*}$ is a fiber bundle, which is known as the Milnor fibration \cite{mil-sing}.
The regular fiber $F = Q^{-1} (1)$ and its transverse intersection $\partial F = Q^{-1} (1) \cap S^{2\ell + 1}_{R}$ with a sphere with large enough radius $R$, are called the \emph{Milnor fiber} and the \emph{Milnor fiber boundary}, respectively.
The Milnor fiber boundary is a closed $(2\ell-1)$-dimensional manifold.

Since each hyperplane is linear, $\cA$ defines a projective hyperplane arrangement $\overline{\cA} = \{\overline{H}_1, \ldots, \overline{H}_n \}$ in $\bC P^{\ell}$. 
We set $U(\cA) = \bC P^{\ell} \setminus \bigcup_{i=1}^{n} \overline{H}_{i}$ the projectivized complement.
It is known that the Hopf map $\pi: M(\cA) \rightarrow U(\cA)$ is a trivial $\bC^{*}$-bundle. Thus, the Betti numbers of the projectivized complement satisfies $b_{k} (M) = b_{k} (U) + b_{k-1} (U)$ for $0 \leq k \leq \ell+1$.

Moreover, its restriction to the Milnor fiber boundary defines a cyclic cover.
 Let $\partial U$ be the boundary of a regular neighborhood of $\bigcup_{i=1}^{n} \overline{H}_{i}$. The manifold $\partial U$ is called the boundary manifold of the projectivized arrangement $\overline{\cA}$.
\begin{proposition} (Lemma 7.5 in \cite{suc-mil})
The restriction $\pi: \partial F \rightarrow \partial U$ is a cyclic $n$-fold cover.
\end{proposition}

Let us consider the case when $\ell = 2$. It is known that $\partial U$ is a closed $3$-manifold obtained from the incidence graph by plumbing construction (see \cite{coh-suc-bound}).
The fundamental group $\pi_{1} (\partial U)$ is generated by meridians $\{ x_i \}$ of each hyperplane and cycles $\{ y_{i} \}$ coming from the incidence graph. 
The characteristic map of the cyclic cover is described as follows (see Proposition 7.6 in \cite{suc-mil}).

\begin{proposition}
The characteristic homomorphism $\omega : \pi_1 (\partial U) \rightarrow \bZ_n$ corrsponding to the $n$-fold cyclic cover $\pi: \partial F \rightarrow \partial U$ is given by $\omega(x_{i})=1$ and $\omega(y_i)=0$.
\end{proposition}

Since $\Hom (\pi_1 (\partial U), \bZ_n) \cong H^{1} (\partial U,\bZ_{n})$, we can consider the characterisitic map belongs to the mod $n$ first cohomology group.
From the argument of Remark 6.7 in \cite{suc-mil}, the characteristic map is expressed as $\omega = i^{*} \omega'$, where $i: \partial U \rightarrow U$ is the inclusion and $\omega' \in H^{1} (U; \bZ_n)$ satisfies $\omega' (x_{i}) = 1$ for each meridian $x_{i}$ of $\overline{H}_i$.

\subsection{The cohomology ring and mod-$2$ resonance in boundary manifolds}\label{sec:resonance}


Next, we review the doubling structure of the cohomology ring of the boundary manifold.
The doubling structure is well-studied in \cite{coh-suc-proj}; however, their arguments imposed restrictions on the coefficient ring (such as characteristic $0$). Therefore, since we wish to work over $\bZ_2$, some care is required.
Nevertheless, most of the arguments run in the same manner as in \cite{coh-suc-proj}.

Let $R$ be commutative ring and $A=\bigoplus_{k=0}^{m}A^{k}$ be a finitely generated graded $R$-algebra. The dual $\overline{A} = \Hom(A, R)$ is also a graded $R$-algebra with the grading $\overline{A}^{k} = \Hom (A^k, R)$.
The dual $\overline{A}$ equips the following $A$-bimodule structure: For $a \in A$, $f \in \overline{A}$, $af, fa \in \overline{A}$ are defined as $(af)(b) = f(ba)$ and $(fa)(b) = f(ab)$. It follows that if $a \in A^{k}$ and $f \in \overline{A}^{\ell}$, then $af, fa \in \overline{A}^{\ell -k}$.
Then, the double $\mathsf{D}(A)$ is defined by the underlying module $\mathsf{D}(A) \cong A \oplus \overline{A}$ with the multiplication $(a,f)\cdot (b,g) = (ab, ag+fb)$.
The double $\mathsf{D}(A)$ has the grading given by $\mathsf{D}(A)^{k} = A^{k} \oplus \overline{A}^{2m-1-k}$.

\begin{theorem} (see also Theorem 4.2 in \cite{coh-suc-proj})
The cohomology ring of the boundary manifold is isomorphic to the double of the cohomology ring of the complement:
\[
H^{*}(\partial U; \bZ_2) \cong \mathsf{D} (H^{*} (U ;\bZ_2)).
\]
\end{theorem}

\proof
Since the statement has been proved only for coefficients in $\bZ$ or $\bC$ in \cite{coh-suc-proj}, we prove that it holds when the coefficient is $\bZ_2$.
By Proposition 2.5 (for $\bZ_2$ coefficient) and Theorem 3.3 in \cite{coh-suc-proj}, it is enough to prove that $H^{*}(U, \partial U; \bZ_2) \subset H^{*}(U; \bZ_2)$ is a square-zero subalgebra.
Since $H^{*}(U; \bZ)$ and $H^{*} (U, \partial U; \bZ)$ are free $\bZ$-modules, the reduction mod $2$ map is surjective. 
By Theorem 4.4 in \cite{coh-suc-proj}, we know that $H^{*} (U, \partial U; \bZ)$ is a square-zero subalgebra of $H^{*} (U; \bZ)$.
The fact that $H^{*}(U, \partial U; \bZ_2)$ is a square-zero subalgebra follows from the commutative diagram below, whose top and bottom rows are exact:
\[
\begin{tikzpicture}
\draw (0,0) node{$H^{*}(U, \partial U; \bZ)$} ;
\draw (3,0) node{$H^{*}(U, \partial U; \bZ_2)$};
\draw (0,-1.5) node{$H^{*}(U; \bZ)$};
\draw (3,-1.5) node{$H^{*}(U; \bZ_2)$};
\draw[->] (0,-0.5) --++(0,-0.5);
\draw[->] (3,-0.5) --++(0,-0.5);
\draw[->] (1.2,0) --++(0.5,0);
\draw[->] (1.0,-1.5) --++(0.8,0);
\draw[->] (4.5,-1.5) --++(0.5,0) node[right]{$0$};
\draw[->] (4.5,0) --++(0.5,0) node[right]{$0$};
\end{tikzpicture}
\]
\endproof




Let $A = H^{*} (U; \bZ_2)$ and $A^{k}=H^{k} (U;\bZ_2)$. From the above proposition, $\mathsf{D}(A) \cong H^{*}(\partial U; \bZ_2)$.
Since $H_{1} (U; \bZ)$ is a free $\bZ$-module, $a^{2} = 0$ for each $a \in A^{1}$ (see Lemma 6.1 in \cite{suc-mod2}). Thus, by fixing an element $a \in A^1$, we have a cochain complex:
\[
(A, a) : 0 \rightarrow A^{0} \xrightarrow{\cdot a} A^{1} \xrightarrow{\cdot a} A^{2} \rightarrow 0.
\]
For each $b \in \overline{A}^{2}$, the element $(a,b) \in \mathsf{D}(A)^{1}$ satisfies $(a,b)\cdot (a,b) = (a^2, ab+ba) = (0,0)$. Thus, we have another cochain complex from the doubling algebra:
\[
(\mathsf{D}(A), (a,b) ) : 0 \rightarrow \mathsf{D}(A)^{0} \xrightarrow{\cdot (a,b) } \mathsf{D}(A)^{1} \xrightarrow{\cdot (a,b)} \mathsf{D}(A)^{2} \xrightarrow{\cdot (a,b)} \mathsf{D}(A)^{3} \rightarrow 0.
\]
Also, the dual $\overline{A} = \Hom (A, \bZ_2)$ defines a cochain complex:
\[
(\overline{A}, a) : 0 \rightarrow \overline{A}^{2} \xrightarrow{\cdot a} \overline{A}^{1} \xrightarrow{\cdot a} \overline{A}^{0} \rightarrow 0.
\] 
A direct computation shows that the following is a short exact sequence of cochain complexes (c.f. Lemma 6.5 and Proposition 6.7 in \cite{coh-suc-proj}):
\[
0 \rightarrow (\overline{A}, a)[-1] \rightarrow (\mathsf{D}(A), (a,b) ) \rightarrow (A, a) \rightarrow 0,
\]
where $(\overline{A},a)[-1]$ is the $(-1)$-shift of the cochain complex. We have the following long exact sequence of the cohomology groups:
\begin{eqnarray*}
0 &\rightarrow& H^{0} (\mathsf{D}(A), (a,b)) \rightarrow H^{0} (A, a) \rightarrow H^{0} (\overline{A}, a) \\
&\rightarrow& H^{1} (\mathsf{D}(A), (a,b)) \rightarrow H^{1} (A, a) \rightarrow H^{1} (\overline{A}, a) \\
&\rightarrow& H^{2} (\mathsf{D}(A), (a,b)) \rightarrow H^{2} (A, a) \rightarrow H^{2} (\overline{A}, a) \\
&\rightarrow& H^{3} (\mathsf{D}(A), (a,b)) \rightarrow 0.
\end{eqnarray*} 
Note that $H^{k} (\overline{A}, a) \cong H^{2-k} (A,a)$. 
If we assume that $a \neq 0$, then each cohomology group satisfies $H^{0}(A, a) = 0$, $H^{1} (A, a) = \bZ_{2}^{d}$ and $H^{2} (A, a) = \bZ_{2}^{\beta+d}$ for some $d$, where $\beta = 1- b_1 (A) + b_{2} (A)$ is the Euler characteristic $\chi(U)$. Thus, the above long exact sequence becomes as follows:
\[
0 \rightarrow \bZ_{2}^{\beta+d} \rightarrow H^{1} (\mathsf{D}(A), (a,b)) \rightarrow \bZ_{2}^{d} \rightarrow \bZ_{2}^{d} \rightarrow  H^{2} (\mathsf{D}(A), (a,b)) \rightarrow \bZ_{2}^{\beta+d} \rightarrow  0.
\]
In particular, by the injectivity of $\bZ_2^{\beta + d} \rightarrow H^{1} (\mathsf{D}(A),(a,b))$, we have the following proposition.
\begin{proposition}\label{prop:resonance}
Suppse that $(a,b) \in A^1 \oplus \overline{A}^2 = \mathsf{D}(A)^1$ satisfies $a \neq 0$. Then, we have 
\[
\dim_{\bZ_2} H^{1}(\mathsf{D}(A), (a,b)) \geq \chi (U).
\]
\end{proposition}

Now, we apply these preparations to the Milnor fiber boundary case.
To apply them, we assume that the number $n$ of hyperplanes is even.
By composing the canonical multiplication $\bZ_n \rightarrow \bZ_2$, we have a surjective homomorphism $\overline{\omega}: \pi_1 (\partial U) \rightarrow \bZ_n \rightarrow \bZ_2$ which belongs to $H^1 (\partial U; \bZ_2)$.
Through the doubling construction, the element $\overline{\omega}$ is expressed as $(\overline{\omega'}, 0)$, where $\overline{\omega'} \in H^{1} (U; \bZ_2)$ satisifies $\overline{\omega'} (x_{i}) = 1$ for each meridian $x_{i}$, particularly, $\overline{\omega'} \neq 0 \in H^{1} (U ; \bZ_2)$. 
Thus, we have the following.
\begin{corollary}\label{cor:resonance}
\[
\dim_{\bZ_2} H^{1} (H^{*} (\partial U; \bZ_2), \overline{\omega}) \geq \chi (U).
\]
\end{corollary}

\section{Torsions in the homology groups of the tower of double coverings} \label{sec:double}
In this section, we review some techniques for studying even torsions in double coverings.

Let $X$ be a connected CW complex and $\omega: \pi_1 (X) \rightarrow \bZ_{2^m}$ be a surjective homomorphism.
Let $p_{\omega}: X^{\omega} \rightarrow X$ be the corresponding $2^{m}$-fold cyclic cover.
By composing $\cdot 2^{m-k}: \bZ_{2^m} \rightarrow \bZ_{2^{k}}$ for $0 \leq k \leq m$, we get a surjective homomorphism $\omega_{k}: \pi_1(X) \rightarrow \bZ_{2^m} \rightarrow \bZ_{2^{k}}$ and corresponding $2^{k}$-fold cyclic cover $p_{\omega_{k}}: X^{\omega_{k}} \rightarrow X$.
We have a tower of double coverings expressed as
\[
X^{\omega} = X^{\omega_m} \rightarrow X^{\omega_{m-1}} \rightarrow \cdots \rightarrow X^{\omega_{1}} \rightarrow X^{\omega_{0}} = X 
\]
The following is proved in \cite{yos-double}.
\begin{proposition}\label{prop:tower} (Lemma 3.4 in \cite{yos-double})
Let $\overline{b}_1(X) = \dim_{\bZ_2} H_{1}(X; \bZ_2)$ be the mod-$2$ Betti number. Then, 
\[
\overline{b}_{1} (X^{\omega_{m}}) \geq \cdots \geq \overline{b}_{1} (X).
\]
\end{proposition}
Let $\eta_{k}: \pi_1 (X^{\omega_{k}}) \rightarrow \bZ_{2}$ be the characteristic map of the double covering $X^{\omega_{k+1}} \rightarrow X^{\omega_{k}}$.
Note that $\eta_{0} = \omega_1$.
From Lemma 3.2 (i) in \cite{yos-double}, $\eta_{k}^2 = 0$ for each $0\leq k \leq m-1$.
Similarly as in Section \ref{sec:resonance}, $\eta_{k}$ defines a cochain complex $(H^{*}(X^{\omega_{k}}; \bZ_2), \eta_{k} )$.
\begin{definition}\label{def:alpha}
We define $\alpha (\omega,k)$ as the first mod-$2$ Aomoto-Betti number of this cochain complex, that is,
\[
\alpha (\omega,k) := \dim_{\bZ_2} H^{1} (H^{*} (X^{\omega_{k}} ; \bZ_2), \eta_{k}).
\]
\end{definition}
Using the Aomoto-Betti number, we can describe the mod-$2$ Betti number of the double cover.

\begin{theorem}\label{thm:mod2betti} (Theorem 3.7 in \cite{yos-double})
The mod-$2$ Betti number of the double cover is described as follows:
\[
\overline{b}_{1} (X^{\omega_{k+1}}) = \overline{b}_{1} (X^{\omega_{k}}) + \alpha (\omega,k).
\]
\end{theorem}

For each double covering $X^{\omega_{k+1}} \rightarrow X^{\omega_{k}}$, the cyclic group $\bZ_2$ acts on $X^{\omega_{k+1}}$ as a decktransformation.
The homology group of the double covering is decomposed into the monodromy eigenspaces:
\[
H_{1} (X^{\omega_{k+1}}; \bC) \cong H_{1} (X^{\omega_{k+1}}; \bC)_{1} \oplus H_{1} (X^{\omega_{k+1}}; \bC)_{-1},
\]
where $H_{1} (X^{\omega_{k+1}}; \bC)_{\pm 1}$ is the $(\pm 1)$-eigenspaces of algebraic monodromy by the decktransformation.
It is known that the $1$-eigenspace is isomorphic to $H_{1} (X^{\omega_{k}}; \bC)$ and the $(-1)$-eigenspace is isomorphic to the corresponding rank one $\bC$-local system homology.
\begin{definition}
Define the dimension of $(-1)$-eigenspace of the homology groups by
\[
\rho(\omega,k) := \dim_{\bC} H_{1}(X^{\omega_{k+1}}; \bC)_{-1}.
\]
\end{definition}
By definition, we have the following equality and inequality;
\begin{eqnarray*}
b_{1} (X^{\omega_{k+1}}) = b_{1} (X^{\omega_{k}}) + \rho (\omega, k), \\
b_1(X^{\omega_{m}}) \geq \cdots \geq b_{1} (X).
\end{eqnarray*}
\begin{definition}
We define the dimension of even torsion summands by
\[
\tau (\omega, k) := \dim_{\bZ_2} (\Tor (H_1 (X^{\omega_k} ; \bZ)) \otimes \bZ_2).
\]
\end{definition}
Ishibashi--Sugawara--Yoshinaga obtained a formula for the number of even torsion summands in the homology group in \cite{isy}, which refines a part of the formula given by Papadima--Suciu (Theorem C in \cite{pap-suc-10}).
\begin{theorem}\label{thm:isy} (Corollary 2.5 in \cite{isy})
The number of even torsion summands satisfies
\[
\tau (\omega, k+1) = \alpha (\omega, k) - \rho (\omega, k).
\]
\end{theorem}

\begin{corollary}\label{cor:torsion}
In the above setting, we additionally suppose that $b_{1} (X^{\omega}) = b_{1} (X)$. Then, 
\begin{enumerate}[(1)]
\item $\alpha(\omega, m-1) \geq \cdots \geq \alpha (\omega, 1) \geq \alpha (\omega, 0). $
\item $\tau (\omega, m) \geq \alpha (\omega, 0)$.
\end{enumerate}
\end{corollary}

\proof
\begin{enumerate}[(1)]
\item 
Since $b_{1} (X^{\omega}) = b_{1} (X)$, we have $b_{1}(X^{\omega_{k}}) = b_{1} (X)$ for each $0 \leq k \leq m$.
By Proposition \ref{prop:tower} and Theorem \ref{thm:mod2betti}, 
\begin{eqnarray*}
b_{1} (X^{\omega_{m-1}}) + \alpha (\omega ,m-1 ) \geq \cdots \geq b_{1} (X) + \alpha (\omega, 0), \\
\alpha(\omega, m-1) \geq \cdots \geq \alpha (\omega, 1) \geq \alpha (\omega, 0). 
\end{eqnarray*}

\item Since $b_{1} (X^{\omega_{k+1}}) = b_{1} (X^{\omega_{k}})$, we have $\rho(\omega,k)=0$ for each $0 \leq k \leq m-1$. Combining this with Theorem \ref{thm:isy} and (1), we have
\[
\tau (\omega, m) = \alpha (\omega, m-1) - \rho (\omega , m-1) \geq \alpha (\omega, 0).
\]
\end{enumerate}
 
\endproof

\section{Proof of the main result.}
We are now ready to prove the main result.
Again, let $\cA = \{H_1, \ldots, H_n\}$ be a central hyperplane arrangement in $\bC^3$.
The assumption on the multiplicity of the intersections is necessary for the following.

\begin{proposition}\label{prop:equal}
Suppose that $(|\cA_X|-2)(\gcd(|\cA_X|, n) -1) =0$ for each $X \in L_2 (\cA)$. 
Then, $b_{1} (\partial F) = b_{1} (\partial U)$.
\end{proposition}

\proof
\begin{equation*}
\begin{aligned}
b_{1} (\partial F) &= \sum_{X \in L_{2} (\cA)} (1 + (|\cA_{X}|-2) \gcd(|\cA_{X}|, n) ) \,   &(\mbox{Theorem \ref{thm:mfb_betti}})  \\
&= \sum_{X \in L_{2} (\cA)} (|\cA_{X}| -1 ) \, &(\mbox{By the assumption}) \\ 
&= \sum_{X \in L_{2} (\cA)} \mu(X) \,& (\mbox{M\"obius function}) \\ 
&= b_{2}(M) \,& (\mbox{Characteristic polynomial}) \\ 
&= b_{1}(U) + b_{2} (U) \,& (\mbox{$M \approx U \times \bC^{*}$}) \\
&= b_{1}(\partial U). \, &(\mbox{Doubling formula})
\end{aligned}
\end{equation*}
\endproof

\begin{remark}
Proposition \ref{prop:equal} also follows from the Alexander polynomial argument (see Theorem 5.2 in \cite{coh-suc-bound}).
\end{remark}

\textit{Proof of Theorem \ref{thm:main}.}
Since $n$ is even, we can assume that $n=s \cdot 2^{m}$ ($m >0$ and $s$ is odd). By multiplying $s \cdot 2^{m-k}$ to the characteristic map, we have a surjective homomorphism 
$\omega_k : \pi_1 (\partial U) \to \bZ_{n} \xrightarrow{\cdot s} \bZ_{2^m} \xrightarrow{\cdot 2^{m-k}} \bZ_{2^k}$. Then, the cyclic cover $\partial F \to \partial U$ decomposes into the following tower of coverings:
\[
\partial F \xrightarrow{s:1} \partial U^{\omega_m} \xrightarrow{2:1} \partial U^{\omega_{m-1}}\xrightarrow{2:1} \cdots \xrightarrow{2:1} \partial U.
\]
Now, let us substitute $X = \partial U$ and $\eta_{0} = \overline{\omega} \in H^{1} (\partial U; \bZ_2)$. Under our assumption, we can apply Corollary \ref{cor:torsion} from Proposition \ref{prop:equal}. By Corollary \ref{cor:resonance}, we have 
\[
\dim_{\bZ_{2}} (\Tor(H_1(\partial U^{\omega_{m}}); \bZ) \otimes \bZ_2)=\tau (\omega, m) \geq
\alpha (\omega, 0) \geq
\chi (U) .
\]
Finally, we apply the transfer argument for the remaining cyclic cover $p: \partial F \to \partial U^{\omega_{m}}$ (see \cite[Section 3.G]{hat} for details). 
The integral transfer satisfies $p_{*} \circ \operatorname{tr} = s \, \mathrm{id}$ in $H_1 (\partial U^{\omega_{m}}; \bZ)$. 
Multiplication by the odd
integer $s$ is an automorphism of the $2$-primary torsion subgroup
of $H_1(\partial U^{\omega_m}; \bZ)$. Hence this subgroup embeds as a direct
summand of the $2$-primary torsion subgroup of
$H_1(\partial F; \bZ)$. In particular, we have 
\[
\dim_{\bZ_{2}} (\Tor(H_1(\partial F); \bZ) \otimes \bZ_2) \geq \dim_{\bZ_{2}} (\Tor(H_1(\partial U^{\omega_{m}}); \bZ) \otimes \bZ_2).
\]
This completes the proof. 
\qed

\section{Examples}
In this section, we give a computation example (see \cite{nem-szi} or Section 5 in \cite{sug-mfb} for the method for obtaining the plumbing graph).
We list all results on the torsion in the first homology group of the Milnor fiber boundary when the number of hyperplanes is eight.

Let $n_k = |\{X \in L_2(\mathcal{A}) \mid |\mathcal{A}_{X}| = k\}| $ denote the number of codimension two intersections of multiplicity $k$. For each tuple $ (n_k)_{k=3}^{8} $, we list the corresponding computational results (Table \ref{tab:8line}). The possible combinatorial types arising from arrangements of projective lines in $ \mathbb{C}P^2 $, and in particular the tuples $ (n_k)_{k=3}^{8} $, have been classified (see \cite{gtv,naz-yos}). Although distinct combinatorial types may share the same tuple $(n_k) $, listing all of them would produce an unmanageably large table. Moreover, in all computations carried out by the author, the homology groups are isomorphic whenever the tuples $ (n_k) $ coincide, even if the combinatorial types differ. Therefore, we summarize the results in a table organized only by the data $ (n_k) $.

In Table \ref{tab:8line}, the assumption $(*)$ stands for the condition 
\[
(|\mathcal{A}_X| - 2)\bigl(\gcd(|\mathcal{A}_X|, n) - 1\bigr) = 0 \,  \mbox{for each $X \in L_2 (\cA)$} \cdots (*).
\]

As the computations show, in every example satisfying assumption $(*)$, the torsion in the homology group coincides with $ \mathbb{Z}_8^{\chi(U)} $. Even in examples not satisfying the assumption, the dimension of the resulting $ \mathbb{Z}_2$-vector space obtained via reduction modulo $2$ equals $ \chi(U)$. Furthermore, in non-satisfying cases, the order of the torsion appears to be related to $ \gcd(|\mathcal{A}_X|, n) $. Understanding the precise relationship between the torsion order and the number $ \gcd(|\mathcal{A}_X|, n) $ remains a topic for future investigation.

\begin{table}[htbp]
  \centering
  \begin{tabular}{|c|c|c|c|c|}
    \hline
$(n_3,n_4,n_5,n_6,n_7,n_8)$ & $\Tor (H_1 (\partial F; \bZ))$ & $\chi (U)$ & Assumption $(*)$ & Remark\\ \hline
$(0,0,0,0,0,1)$  & $0$ & $-6$ && Pencil case \\ \hline
$(0,0,0,0,1,0)$ & $0$ & $0$ &$\circ$& Near-pencil case\\ \hline
$(1,0,0,1,0,0)$ & $\bZ_{4}^{4} $ & $4$ & &  \\ \hline
$(0,0,0,1,0,0)$ & $\bZ_{4}^{4} \oplus \bZ_{8}$ & $5$ &&  \\ \hline
$(6,1,0,0,0,0)$ & $\bZ_{2}^{2} \oplus \bZ_{8}^{4}$ & $6$& &  \\ \hline
$(3,2,0,0,0,0)$ & $\bZ_{2}^{4} \oplus \bZ_{8}^{2}$ & $6$& &  \\ \hline
$(3,0,1,0,0,0)$ & $\bZ_{8}^{6}$ & $6$ & $\circ$ &  \\ \hline
$(0,1,1,0,0,0)$ & $\bZ_{2}^{2} \oplus \bZ_{8}^{4}$ & $6$& &  \\ \hline
$(8,0,0,0,0,0)$ & $\bZ_{8}^{7}$ & $7$ & $\circ$ & MacLane arrangement \\ \hline
$(5,1,0,0,0,0)$ & $\bZ_{2}^{2} \oplus \bZ_{8}^{5}$ & $7$ &&  \\ \hline
$(2,2,0,0,0,0)$ & $\bZ_{2}^{4} \oplus \bZ_{8}^{3}$ & $7$ &&  \\ \hline
$(2,0,1,0,0,0)$ & $\bZ_{8}^{7}$ & $7$ &$\circ$&  \\ \hline
$(7,0,0,0,0,0)$ & $\bZ_{8}^{8}$ & $8$ & $\circ$ &  \\ \hline
$(4,1,0,0,0,0)$ & $\bZ_{2}^{2} \oplus \bZ_{8}^{6}$ & $8$ &&  \\ \hline
$(1,2,0,0,0,0)$ & $\bZ_{2}^{4} \oplus \bZ_{8}^{4}$ & $8$ &&  \\ \hline
$(1,0,1,0,0,0)$ & $\bZ_{8}^{8}$ & $8$ &$\circ$&  \\ \hline
$(6,0,0,0,0,0)$ & $\bZ_{8}^{9}$ & $9$ &$\circ$&  \\ \hline
$(3,1,0,0,0,0)$ & $\bZ_{2}^{2} \oplus\bZ_{8}^{7}$ & $9$ &&  \\ \hline
$(0,2,0,0,0,0)$ & $\bZ_{2}^{4} \oplus\bZ_{8}^{5}$ & $9$ &&  \\ \hline
$(0,0,1,0,0,0)$ & $\bZ_{8}^{9}$ & $9$ &$\circ$&  \\ \hline
$(5,0,0,0,0,0)$ & $\bZ_{8}^{10}$ & $10$ &$\circ$&  \\ \hline
$(2,1,0,0,0,0)$ & $\bZ_2^{2} \oplus \bZ_{8}^{8}$ & $10$ &&  \\ \hline
$(4,0,0,0,0,0)$ & $\bZ_{8}^{11}$ & $11$ &$\circ$&  \\ \hline
$(1,1,0,0,0,0)$ & $\bZ_2^{2} \oplus \bZ_{8}^{9}$ & $11$ &&  \\ \hline
$(3,0,0,0,0,0)$ & $\bZ_{8}^{12}$ & $12$ &$\circ$&  \\ \hline
$(0,1,0,0,0,0)$ & $\bZ_2^{2} \oplus \bZ_{8}^{10}$ & $12$ &&  \\ \hline
$(2,0,0,0,0,0)$ & $\bZ_{8}^{13}$ & $13$ &$\circ$&  \\ \hline
$(1,0,0,0,0,0)$ & $\bZ_{8}^{14}$ & $14$ &$\circ$&  \\ \hline
$(0,0,0,0,0,0)$ & $\bZ_{8}^{15}$ & $15$ &$\circ$&Generic arrangement  \\ \hline

  \end{tabular}
  \caption{Computation of $\Tor (H_{1} (\partial F; \bZ))$ for $n=8$.}
  \label{tab:8line}
\end{table}


\begin{thebibliography}{99}


\bibitem[CS06]{coh-suc-proj}
D. C. Cohen, A. Suciu, 
\emph{Boundary manifolds of projective hypersurfaces},
Adv. Math,
\textbf{206}, no. 2, 538--566, (2006).

\bibitem[CS07]{coh-suc-bound}
D. C. Cohen, A. Suciu, 
\emph{The boundary manifold of a complex line arrangement},
Groups, homotopy and configuration spaces, 105--146.
Geom. Topol. Monogr., 13
Geometry \& Topology Publications, Coventry, 2008


\bibitem[Di16]{dim-hyp}
A. Dimca, 
\emph{Hyperplane arrangements},
Universitext, Springer, Cham, 2017, xii+200 pp.






\bibitem[GTV03]{gtv}
D. Garber, M. Teicher, U. Vishne,
\emph{$\pi_1$-classification of real line arrangements with up to eight lines},
Topology, \textbf{42} (1), 2003, 265-289.



\bibitem[Ha02]{hat}
A.~Hatcher, 
\emph{Algebraic Topology},
Cambridge University Press, Cambridge, 2002, xii+544 pp.




\bibitem[ISY25]{isy}
S. Ishibashi, S. Sugawara, M. Yoshinaga,
\emph{Betti numbers and torsions in homology groups of double coverings},
Adv. in Appl. Math., \textbf{162} (2025), 102790.



\bibitem[MP03]{mic-pic}
F. Michel, A. Pichon,
\emph{On the boundary of the Milnor fibre of nonisolated singularities},
Int. Math. Res. Not. \textbf{43}, 2305-2311, (2003).

\bibitem[Mi68]{mil-sing}
J. Milnor,
\emph{Singular points of complex hypersurfaces},
Ann. Math. Studies, 
vol. \textbf{61}, Princeton University Press, Princeton (1968).


\bibitem[NY10]{naz-yos}
S. Nazir, M. Yoshinaga,
\emph{On the connectivity of the realization spaces of line arrangements},
Ann. Sc. Norm. Super. Pisa Cl. Sci. (5),
Vol. XI (2012), 921-937.

\bibitem[NS12]{nem-szi}
A. N\'emethi, A. Szil\'ard, 
\emph{Milnor fiber boundary of a non-isolated surface singularity}, 
Lecture Notes in Math. vol. 2037, 
Springer-Verlag, Berlin Heidelberg (2012)



\bibitem[OS80]{orl-sol}
P. Orlik, L. Solomon, 
\emph{Combinatorics and topology of complements of hyperplanes},
Invent. Math., {\bf 56} (1980), no. 2, 167--189.

\bibitem[OT92] {orl-ter}
P. Orlik, H. Terao, 
\emph{Arrangements of Hyperplanes}, 
Grundlehren Math. Wiss.  300, Springer-Verlag New York, 1992. 


\bibitem[PS10]{pap-suc-10}
S. Papadima, A. Suciu,
\emph{The spectral sequence of an equivariant chain complex and homology with
local coefficients}, 
Trans. Amer. Math. Soc., {\bf 362} (5),
2685--2721.

\bibitem[PS17]{pap-suc-17}
S. Papadima, A. Suciu,
\emph{The Milnor fibration of a hyperplane arrangement: from modular resonance to algebraic monodromy}, 
Proc. London Math. Soc., {\bf 114} (2017),
no. 6, 961--1004.



\bibitem[Suc14] {suc-mil}
A. Suciu, 
\emph{Hyperplane arrangements and Milnor fibrations}, 
Ann. Fac. Sci. Toulouse Math. (6)
{\bf 23} (2014) no. 2, 417--481.



\bibitem[Suc23]{suc-mod2}
A. Suciu,
\emph{Cohomology, Bocksteins, and resonance varieties in characteristic 2},
Compactification, configurations, and cohomology, 131-157.
Contemp. Math., 790,
Amer. Math. Soc, Providence, RI, 2023.

\bibitem[Sug25]{sug-mfb}
S. Sugawara,
\emph{First homology groups of the Milnor fiber boundary for generic hyperplane arrangements in $\bC^{3}$}, arXiv:2404.01555, to appear in Bull. London Math. Soc.



\bibitem[Yo20]{yos-double}
M. Yoshinaga, 
\emph{Double coverings of arrangement complements and 2-torsion in Milnor fiber homology}, 
Eur. J. Math. \textbf{6} (2020), no. 3, 1097-1109.



\end{thebibliography}
\end{document}